\documentclass[12pt]{amsart}

\usepackage{amssymb,latexsym, mathrsfs, pifont}

\newtheorem{theorem}{Theorem}

\newtheorem{lemma}{Lemma}

\newtheorem{defn}{Definition}

\newtheorem{proposition}[theorem]{Proposition}

\newtheorem{corollary}[theorem]{Corollary}

\def\RR{{\mathbb R}}

\def\QQ{{\mathbb Q}}

\def\cL{{\mathscr L}}

\def\cH{{\mathcal H}}

\def\slr{$SL(2,\mathbb{R})$}

\begin{document}

\title{On unipotent flows in $\mathcal{H}(1,1)$}

\author{Kariane Calta and Kevin Wortman}

\begin{abstract}

We study the action of the horocycle flow on the moduli space of abelian differentials in genus two. In particular, we exhibit a classification of a specific class of probability measures that are invariant and ergodic under the horocycle flow on the stratum $\cH(1,1)$.

\end{abstract}

\maketitle

\begin{section}{Introduction}\label{S:intro}

We begin by briefly describing translation surfaces in genus two and providing examples of surfaces which are contained in the support of the measures of interest. Then we recall from \cite{calta} the description of the submanifolds on which our measures are supported. An alternate description is given by McMullen \cite{mc1}. For a detailed exposition of translation surfaces and a historical survey of salient results, we refer the reader to Masur-Tabachnikov \cite{mt} or Masur-Hubert-Schmidt-Zorich \cite{hsz}. There are numerous sources of reference for Ratner's measure classification theorem for homogeneous spaces; some of these include the following works of Ratner, Margulis-Tomanov, and Morris: \cite{r1}, \cite{r2}, \cite{r3}, \cite{r4}, \cite{M-T}, and \cite{Mo}.

\begin{subsection}{Closed submanifolds in $\cH(1,1)$}

The moduli space of translation surfaces of genus $g$ is stratified according to the number and order of conical singularities of the surfaces. The Riemann-Roch theorem implies that if an abelian differential on a genus $g$ surface has zeros of orders $m_1, \cdots, m_k$, then $\lbrace m_1, \cdots, m_k \rbrace$ is a partition of $2g-2$. That is, $\sum_{i=1}^k m_i =2g-2$.  If we let $\alpha=\lbrace m_i \rbrace_{i=1}^k$, then we have

\begin{defn} $\cH(\alpha)$ is the stratum of the moduli space of genus $g$ translation surfaces consisting of those surfaces with abelian differentials having $k$ zeros of orders $m_i$. \end{defn}

\medskip

{\bf Example.} In genus two, there are two strata, $\cH(2)$ and $\cH(1,1)$. Surfaces in $\cH(2)$ have one singularity of order 2 or total angle $ 6 \pi$. An example of a surface in $\cH(2)$ is the regular octagon with parallel sides identified. Other examples include the L-tables described by Calta and McMullen. (\cite{calta}, \cite{mc1}) On the other hand, surfaces in $\cH(1,1)$ have two simple singularities, each of total angle $4 \pi$. The regular decagon with parallel sides identified is in this stratum, as are the Z-tables described by Calta and McMullen. (\cite{calta}, \cite{mc4}) Our classification theorems pertain to the stratum $\cH(1,1)$.

\medskip

The group $SL(2,\RR)$ acts on the moduli space of translation
surfaces and preserves each stratum $\cH(\alpha)$.

\begin{defn} A translation surface $S$ is said to be a {\em lattice surface} if $\text{{\em Stab}}_{SL(2,\RR)}(S)$ is a lattice in $SL(2,\RR)$.

\end{defn}

{\bf Example.} A theorem of Veech \cite{veech} shows that the regular decagon is a lattice surface. (In fact, the same theorem implies that any regular $2n$-gon with parallel sides identified is a lattice surface.)

We now recall a definition from \cite{calta} and describe the submanifolds of interest.

A direction $v$ on a translation surface is said to {\em completely periodic} if in the direction $v$, the surface decomposes as a union of cylinders of closed, parallel trajectories bounded by saddle connections.  We say that a surface is completely periodic if any direction in which there exists a closed saddle connection is completely periodic. A theorem of Veech \cite{veech} implies that any lattice surface is completely periodic, although the converse is false and in fact there are counterexamples in $\cH(1,1)$. (See \cite{calta} for a concrete examples.)

In genus two, the holonomy field of a completely periodic surface is quadratic \cite{calta}. We will use this fact in the description of our submanifolds.

Given a surface $S \in \mathcal{H}(1,1)$ and a cylinder decomposition of $S$, let $w_i, h_i$ and $t_i$ denote the width, height, and twist parameters of the cylinders. Number the cylinders so that $w_3 = w_2 + w_1$ and define $s_i=h_i + h_3$ and $\tau_i=t_i + t_3$ for $i=1,2$. (See \cite{calta} for details.) We have the following result:

\begin{theorem} Let $\mathscr{L}_m$ be the set of completely periodic surfaces in $\cH(1,1)$ that can be rescaled so that each surface $S$ has a cylinder decomposition with parameters $w_i, s_i, \tau_i \in \QQ(\sqrt{d})$ that satisfy:

\begin{align} w_1 \bar{s_1}&+ w_2 \bar{s_2}=0 \notag \\
w_1 \bar{\tau_1 } &+ w_2 \bar{\tau_2}=0 \notag \\
w_1 s_1 &+ w_2 s_2 = m. \notag \end{align}

(Here, $\bar{}$ denotes conjugation in $\QQ(\sqrt{d})$ and $m \in \QQ(\sqrt{d})$. Note that $m$ is the area of $S$.)

Then $\mathscr{L}_m$ is a closed, $SL(2,\RR)$-invariant submanifold of $\mathcal{H}(1,1)$. Any primitive lattice surface in $\cH(1,1)$ is contained in some $\mathscr{L}_m$.

\end{theorem}

{\em Note.} In the language of McMullen \cite{mc1}, the $\mathscr{L}_m$ are spaces of eigenforms of fixed discriminant.

\medskip

Since the regular decagon is a lattice surface, it is contained in a submanifold $\mathscr{L}_m$.  It should be noted here as well that McMullen showed in \cite{mc3} that the regular decagon is the only primitive lattice surface in $\cH(1,1)$.

In addition to the action of $SL(2, \RR)$, there is a local pseudo-action of $\RR^2$ on each submanifold $\mathscr{L}_m$. Given a surface $S \in \mathscr{L}_m$, we can define a new surface $vS$ for sufficiently small vectors $v \in \RR^2$. $S$ can be realized as a union of polygons in $\RR^2$ glued along parallel sides; under this identification, there are two equivalence classes of vertices of the polygons, $[p]$ and $[q]$, which give rise to the two singularities of $S$. Choose an equivalence class $[p]$ and add to each vertex in $[p]$ the same vector $v$. Then $vS$ is the translation surface realized as the union of polygons with vertices $p+v$ and $q$, and edges connecting $p+v$ and $q$ for every pair $p$ and $q$ which were connected by an edge in the polygonal decomposition of $S$. (See \cite{calta} for details.) This construction changes the relative homology of $S$ while fixing the absolute homology.

Of particular interest for us in this paper will be translating singularities in $\mathscr{L}_m$ by elements of the group of horizontal vectors $X=\lbrace (x,0) \in \RR^2 \rbrace$. Any surface $S \in \mathscr{L}_m$ is associated with a maximal open interval $I_S \subseteq X$ such that $xS$ is a well defined surface in $\mathcal{H}(1,1)$ -- and hence in $\mathscr{L}_m$ -- for any $x \in I_S$. We let $\overline{I_S} \subseteq X$ be the closure of $I_S$, and we call any vector $x \in \overline{I_S} \setminus I_S$ a \emph{horizontal saddle connection} on $S$.

For example, it is a straightforward task to construct a Z-table $S \in \mathscr{L}_m$ such that $wS \in \mathscr{L}_m$ for all $w \in X$ with $\left\| w \right\| < \ell $ for some fixed $\ell >0$, but such that $x \in X$ is a horizontal saddle connection between distinct singularities on $S$ with $\left\| x \right\|=\ell$. Thus, 
$x S \notin \mathcal{H}(1,1)$ as $x$ has the effect of identifying the two distinct singularities on $S$. 

With $S$ and $x$ as in the above paragraph, $x$ is a horizontal saddle connection on $S$.
Note that through this construction, we could have that $xS$ is an L-table and thus is contained in $\mathcal{H}(2)$. Alternatively, $xS$ could be a table made of two squares identified on a single vertex and such that the sides of each individual square are identified as to form a torus. That is, $xS$ is two tori joined at a single point.

{\em Remark.} Although the group $G=SL(2,\RR)$ acts on each $\cL_m$, there is only a local pseudo-action by $\RR^2$, and so we do not obtain an action by the group $G \ltimes \RR^2$. However, for any point $S \in \mathscr{L}_m$, there is a neighborhood of the identity $\mathcal{O} \subseteq G \ltimes \RR^2$ such that the map $(g,v) \mapsto v(gS)$ defined on $\mathcal{O}$ is a homeomorphism. Thus, there is a local pseudo-action of $G \ltimes \RR^2$ on $\mathscr{L}_m$.

\end{subsection}

\begin{subsection}{Statement of the main result.}

Let $U\leq SL(2,\RR)$ be the subgroup of upper-triangular
unipotent matrices.

\begin{defn} A measure $\mu$ on $\mathscr{L}_m$ is \emph{horizontally
invariant} if the space
$$\mathscr{L}_m^X =\{\,M\in \mathscr{L}_m \mid xM \in \mathcal{H}(1,1) \text{ for all }x\in X\,\}$$ is
conull, and if $\mu$ is $X$-invariant. \end{defn}

\begin{theorem}\label{T:main}

Let $\mu$ be an ergodic $U$-invariant Borel probability measure
supported on some $\mathscr{L}_m$. Further, assume that either
$\mu$ is invariant under the diagonal subgroup of $SL(2,\RR)$, or that $\mu$ is not horizontally invariant, or that $\mu$-almost every $M \in \mathscr{L}_m$ contains no horizontal saddle connections. Then either
\begin{quote}
(i) $\mu$ is the unique
$SL(2,\RR)$-invariant, ergodic probability measure
with $\text{Supp}(\mu)=\mathscr{L}_m$,
\newline (ii) $\mu$ is
arclength on the $U$-orbit of a $U$-periodic surface,
\newline (iii) there is an $x \in X$ such that $x_*\mu = \nu_{10}$
where $\nu _{10}$ is Haar measure on the
closed $SL(2,\RR)$-orbit of the regular decagon, or
\newline (iv) there is an $x \in X$ such that $ x$ is a horizontal saddle connection for every surface in the support of $\mu$. Furthermore, $\mu$ is invariant under the natural action of $x^{-1}SL(2,\RR)x$
on the support of $\mu$.
\end{quote}

\end{theorem}

{\em Remark.} If in case (iv) there is some $S \in \text{Supp}(\mu)$ with $xS \in \mathcal{H}(2)$, then $x _*\mu$ is an ergodic, $SL(2,\mathbb{R})$-invariant measure on $\mathcal{H}(2)$. It follows from McMullen's Theorem 1.5 of \cite{mc2} that $x_*\mu$ is Haar measure on a lattice surface in $\mathcal{H}(2)$.

If there is some $S \in \text{Supp}(\mu)$ such that $xS$ is two tori joined at a point, then $x _*\mu$ is an ergodic, $SL(2,\mathbb{R})$-invariant measure on a product of two moduli spaces of a torus. As this space is homogeneous, Ratner's measure classification theorem applies.

\end{subsection}

\begin{subsection}{Basis of the proof.}

The proof utilizes a framework of results that includes an analysis of the structure of moduli space in genus two, the ideas involved in the proof of Ratner's measure classification theorem, and the proof techniques of an analog of Ratner's theorem for spaces of branched covers of lattice surfaces given by Eskin-Marklof-Morris \cite{E-M-M}.

  Many of the techniques and ideas used in the proof of Theorem~\ref{T:main} are similar to those used in the proof of Ratner's theorems, and to those used in the Eskin-Marklof-Morris measure classification theorem for unipotent flows on spaces of branched covers of lattice surfaces. Examples of these techniques and ideas include the notion of polynomial divergence in unipotent flows, measurement of the directions of transverse divergence under unipotent flows, and entropy inequalities for actions of semisimple elements. In fact, the loose outline for the proof of our theorem is modeled on the proof
of Ratner's theorem provided by Margulis-Tomanov \cite{M-T}, as in \cite{E-M-M}.

 Furthermore, in independent work, McMullen and Calta located and explicitly described the submanifolds $\mathscr{L}_m$ (\cite{calta}, \cite{mc1}) and McMullen classified the ergodic, $SL(2,\RR)$-invariant ergodic probability measures on the moduli space in genus two \cite{mc2}. We make essential use of this classification as well.

\smallskip \noindent \textbf{Acknowledgements.} We thank Alex
Eskin for showing this problem to us and for explaining the argument used in Section~\ref{S:UX}.

We thank John Smillie for his substantial help,
Matthew Bainbridge for pointing out an error in an earlier draft, and Alireza Salehi Golsefidy, Amir Mohammadi, and Barak Weiss for helpful conversations.

\end{subsection}

\end{section}

\section{Horizontal saddle connections and the support of $\mu$}\label{S:general}

We fix an ergodic $U$-invariant probability measure on some
$\mathscr{L}_m$, and denote it by $\mu$.

\subsection{Well-defined horizontal translations}\label{s:X}

Recall that $X \leq \mathbb{R}^2$ is the subgroup of $\RR^2$ of horizontal
vectors. It will often be convenient for us to identify a vector
in $X$ with the real value of its nontrivial coordinate.

In what follows, for any $M \in \mathscr{L}_m$ and any interval $Z
\subseteq X$ with $vM \in \mathscr{L}_m$ for all $v \in Z$, we denote $\cup_{v \in Z}vM$ as $ZM$.

In particular, for any submanifold $\mathscr{L}_m$ and any $s>0$, we define $$\text{HC}(s)=\{\,M\in \mathscr{L}_m \mid
(-s,s)M\subseteq \mathscr{L}_m \mbox{ and } [-s,s]M \nsubseteq
\mathcal{H}(1,1) \,\}$$
where a real number $t \in (-s,s)$ is identified with $(t,0) \in
X$, as described above.

Note that since $U$ and $X$ commute, each $\text{HC}(s)$ is
$U$-invariant, and that
$$\mathscr{L}_m=\mathscr{L}_m^X  \coprod \Big(\coprod_{s >0}\text{HC}(s)\Big)$$

Since $\mathscr{L}_m^X$ is also $U$-invariant, it follows from the
ergodicity of $\mu$ that either $\mathscr{L}_m^X$ is conull or $\text{Supp}(\mu) \subseteq
\text{HC}(s)$ for some $s>0$.

\begin{lemma}\label{l:X_c} There is an open interval $X_0
\subseteq X$, and a conull set $\mathscr{C}_m \subseteq
\mathscr{L}_m$ such that $X_0 M \subseteq \mathscr{L}_m$ for any
$M \in \mathscr{C}_m$. Furthermore, we may assume that $X_0 = -X_0$.
\end{lemma}

\begin{proof}

If $\mathscr{L}_m^X$ is conull, let $X_0=X$ and
$\mathscr{C}_m=\mathscr{L}_m^X$. On the other hand, if $\text{Supp}(\mu)
\subseteq \text{HC}(s)$, let $X_0 = (-s,s)$ and $\mathscr{C}_m=\text{HC}(s)$.
\end{proof}

\section{Transverse divergence}\label{S:transverse}

We denote the set of vertical vectors in the plane as:
$$Y=\{\,(0,y) \in \mathbb{R} ^2 \mid y\in \mathbb{R} \,\}.$$
It will be convenient to write $\RR^2$ as $XY$ to emphasize the
difference between the action of a vertical vector and a
horizontal vector. The difference between $X$ and $Y$ is
significant because although $X$ is $U$-invariant, $Y$ is not.

 Now $SL(2,\RR)$ acts on $XY$ so that we may form the semidirect product $SL(2,\RR) \ltimes XY$, a $5$-dimensional Lie group.
We embed $SL(2,\RR) \ltimes \RR^2$ into $SL(3,\RR)$:
$$\begin{bmatrix} * & * & * \\ * & * & * \\ 0 & 0 & 1 \end{bmatrix}$$ We multiply elements of $SL(2,\RR) \ltimes \RR^2$ by multiplying their images in $SL(3,\RR)$.

For any $M \in \mathscr{L}_m$, there is an open
neighborhood of the identity $\mathcal{O}_M \subseteq
SL(2,\RR)\ltimes XY$ such that the map
$$(g,xy) \mapsto xy(g M)$$ defines a homeomorphism of $\mathcal{O}_M$
onto its image in $\mathscr{L}_m$.

For any $s \in \mathbb{R}$, we let $x^s=(s,0) \in \mathbb{R}^2$
and \[ a^s = \begin{pmatrix}
e^s & 0  \\
0 & e^{-s}
\end{pmatrix}
\]

We let $B \leq SL(2,\RR)$ be the subgroup of upper-triangular matrices
and we define $\text{Stab}_{B{X_0}}(\mu
)$  to be the set of all $h \in BX_0$ such that $hM$ is well-defined for almost every $M \in \mathscr{L}_m$ and such that $  h_*\mu =\mu $. Note that
$\text{Stab}_{B{X_0}}(\mu )$ is merely a set, and may not be a
group if $X_0 \neq X$.

The following lemma is a key tool used in the proof for
Theorem~\ref{T:main}. It is essentially Proposition 4.4 from
\cite{E-M-M}.

First, we establish some notation. 
Let 

 $$g_k=\begin{bmatrix} a_k & b_k \\ c_k & d_k \\ \end{bmatrix}$$
and 

$$f(g_k,t_k)= \frac{t_k a_k - b_k}{d_k - t_k c_k}.$$

\begin{lemma}\label{l:transverse}

Given $\varepsilon > 0$, there exists a compact set $\Omega _\varepsilon \subseteq
\mathscr{L}_m$ with $\mu (\Omega _\varepsilon ) >1- \varepsilon$
and some $0<\delta<1$ for which $x^\delta \in X_0$, such that if

\begin{quote} (i) $\{M_k\}$ is a convergent sequence in
$\Omega _\varepsilon$ with $\lim
_{k \to \infty} M_k=M'$,\\
(ii) $(g_k,x_ky_k)M_k=M'$ for $(g_k,x_ky_k) \in
SL(2,\RR)\ltimes XY$ with $\lim
_{k \to \infty} (g_k,x_ky_k)=1$,\\
(iii) either $|c_k| \leq |y_k|$ for all $k$ or $|y_k| \leq |c_k|$ for all $k$, \\
(iv) $c_ky_k\neq 0$ for all $k$, and \\
(v) $t_k = \min \{\, \frac{\delta}{|c_k|} \,,\,
\frac{\delta}{ |y_k| } \,\}$,
\end{quote}

\noindent then $\emph{Stab}_{BX_0}(\mu )^\circ$ contains an
element of $AX_0$. More precisely, if $|y_k|\leq|c_k|$ for all $k$ then $\emph{Stab}_{BX_0}(\mu )^\circ$ contains
an element of $$a^{-\log (1\pm \delta)}x^{[\frac{-\delta}{1-\delta},\frac{\delta}{1-\delta}]}$$ and if $|c_k|\leq|y_k|$ for all $k$ then $\emph{Stab}_{BX_0}(\mu )^\circ$ contains
an element of $$a^{[-\log (1 + \delta),-\log(1-\delta)]}\big[x^{[\frac{-\delta}{1-\delta},\frac{-\delta}{1+\delta}]} \cup
x^{[\frac{\delta}{1+\delta},\frac{\delta}{1-\delta}]}\big] $$

\end{lemma}

\begin{proof} As the proof is essentially that of \cite{E-M-M}
Proposition 4.4, we will only sketch the argument and refer to
\cite{E-M-M} for the details.

Let $\Omega _\epsilon$ be a compact uniformly generic set for the
action of $U$ on $\mathscr{L}_m^X$. There is a compact set $K
\subseteq \mathscr{L}_m$ with measure arbitrarily close to $1$ such that $K \cap (bx)K =\emptyset$ for any $bx \in
BX_0-\text{Stab}_{BX_0}(\mu)$.

Let \[ u^t = \begin{pmatrix}
1 & t  \\
0 & 1
\end{pmatrix}
\]

Since $M' \in \Omega _\varepsilon$,
\begin{align}\label{e:1}d(u^{t}M', K)<\varepsilon\end{align} for most
$t \in [-t_k,t_k]$.

Let $f(g,t)=\frac{ta-b}{d-tc}$. Since the derivative of $f$ is
bounded on compact intervals, it follows that
\begin{align}\label{e:2}d(u^{f(g_k,t)}M_k, K)<\varepsilon\end{align} for most $t \in
[-t_k,t_k]$ as well.

Recall that $\mathcal{O}_{M'}$ is an open neighborhood of the
identity in $SL(2,\RR)XY$ which is homeomorphic to its image in
$\mathscr{L}_m$ via the map $(g,xy) \to xygM'$. Now, there is also an
open neighborhood, $\mathcal{O}$, of $M'$ in $\mathcal{O}_{M'}M'$
and some $\delta >0$ such that for any $N \in
\mathcal{O}$ and any $v \in XY$ with $||v||<\delta$ we have $vN
\in \mathcal{O}_{M'}M'$. We may assume $M_k \in \mathcal{O}$ for
all $k$.

We claim that after passing to a subsequence, $u^{f(g_k,t_k)}g_kx_ky_ku^{-t_k}$
converges to the element of $\text{Stab}_{BX_0}(\mu )^\circ$ as described in the conclusion of the statement of this lemma.
Call this limit
$\beta_\delta \in BX_0$.


Now we prove the claim.
Our goal is to compute $u^{f(g_k,t_k)} g_k x_k y_k u^{-t_k}$. 

Since we assumed that $g_kx_ky_k \to 1$, we know that $a_k, d_k \to 1$ and $c_k, b_k, x_k, y_k \to 0$. 

For convenience, we let  $f(g_k,t_k)= f_k$.

First note that \begin{align*}g_kx_ky_k & = \begin{bmatrix}
 a_k & b_k & 0 \\
 c_k & d_k & 0 \\
  0 & 0 & 1 \end{bmatrix}
  \begin{bmatrix}
 1 & 0 & x_k \\
 0 & 1 & 0 \\
  0 & 0 & 1 \end{bmatrix}
  \begin{bmatrix} 1 & 0 & 0 \\
 0 & 1 & y_k \\
  0 & 0 & 1 \end{bmatrix}
  \\
  &=\begin{bmatrix}
 a_k & b_k & a_kx_k + b_ky_k \\
 c_k & d_k  & c_kx_k+d_ky_k \\
  0 & 0 & 1 \end{bmatrix}
  \end{align*}

Now note that \begin{align*}g_kx_ky_k u^{-t_k} & =\begin{bmatrix}
 a_k & b_k & a_kx_k + b_ky_k \\
 c_k & d_k  & c_kx_k+d_ky_k \\
  0 & 0 & 1 \end{bmatrix}
  \begin{bmatrix}
 1 & -t_k & 0 \\
 0 & 1 & 0 \\
  0 & 0 & 1 \end{bmatrix}
  \\
& =\begin{bmatrix}
 a_k & b_k-t_ka_k & a_kx_k + b_ky_k \\
 c_k & d_k-t_kc_k  & c_kx_k+d_ky_k \\
  0 & 0 & 1 \end{bmatrix}
 \end{align*}

Next, note that $$ u^{f(g_k,t_k)} g_k x_k y_k u^{-t_k}  =  \begin{bmatrix}
 1 & f_k & 0 \\
 0 & 1 & 0 \\
 0 & 0 & 1 \end{bmatrix}
\begin{bmatrix}
 a_k & b_k-t_ka_k & a_kx_k + b_ky_k \\
 c_k & d_k-t_kc_k  & c_kx_k+d_ky_k \\
 0 & 0 & 1 \end{bmatrix} $$

  $$ = \begin{bmatrix}
 a_k+c_kf_k & (b_k-t_ka_k)+f_k(d_k-t_kc_k) & (a_kx_k + b_ky_k)+f_k(c_kx_k+d_ky_k) \\
 c_k & d_k-t_kc_k  & c_kx_k+d_ky_k \\
 0 & 0 & 1 \end{bmatrix} $$

  $$ = \begin{bmatrix}
 a_k+c_kf_k & 0 & (a_kx_k + b_ky_k)+f_k(c_kx_k+d_ky_k) \\
 c_k & d_k-t_kc_k  & c_kx_k+d_ky_k \\
 0 & 0 & 1 \end{bmatrix} $$

We have by definition that $t_k=\min \lbrace \frac{\delta}{|c_k|}, \frac{\delta}{|y_k|} \rbrace$.

Suppose that $|y_k| \leq |c_k|$. Then $t_k=\frac{\delta}{|c_k|}$ so $t_kc_k=\pm\delta$. 

Thus, $d_k-t_kc_k=d_k \pm \delta$ and $d_k \pm \delta \to 1 \pm \delta$. So after passing to a subsequence,  $ a_k+c_kf_k \to 1/(1 \pm \delta)$ since our matrix has determinate equal to $1$. 

Since $c_k \to 0$, after projecting $\pi :SL(2,\RR)XY \rightarrow SL(2,\RR)$ we have $\pi(u^{f(g_k,t_k)}g_kx_ky_ku^{-t_k}) \to a^{-\log(1\pm \delta)}$.

As for the ``vector-part'' of the limit, note that $c_kx_k+d_ky_k \to 0$, so we are only left to determine the limit of 
$(a_kx_k + b_ky_k)+f_k(c_kx_k+d_ky_k)$

Since $x_k \to 0$ and $y_k \to 0$,

\begin{align*}\lim[(a_kx_k + b_ky_k)+f_k(c_kx_k+d_ky_k)] & =\lim[f_k(c_kx_k+d_ky_k)] \\
& = \lim[\frac{(t_ka_k-b_k)(c_kx_k+d_ky_k)}{d_k-t_kc_k}] \\
& = \frac{\lim[(t_ka_k-b_k)(c_kx_k+d_ky_k)]}{\lim[d_k-t_kc_k]}
\end{align*}

Since $d_k \to 1$ and $t_k=\delta /|c_k|$, we have $\lim[d_k-t_kc_k ]= 1\pm\delta$.

Furthermore, \begin{align*} \lim[(t_ka_k-b_k)(c_kx_k+d_ky_k)] &= \lim[t_ka_kc_kx_k + t_ka_kd_ky_k-b_kc_kx_k-b_kd_ky_k] \\
& = \lim[\pm \delta a_kx_k + a_kd_k\delta \frac{y_k}{|c_k|}-b_kc_kx_k-b_kd_ky_k] \\
& = \lim[0 + \delta \frac{y_k}{|c_k|}-0-0] \\
& = \delta \lim \frac{y_k}{|c_k|} \\
\end{align*}

Since $|y_k| \leq |c_k|$, the above limit is contained in the interval $[-\delta, \delta]$. Thus,

$$\lim[(a_kx_k + b_ky_k)+f_k(c_kx_k+d_ky_k)] $$ is contained in the interval $$[\frac{-\delta}{1-\delta},\frac{\delta}{1-\delta}]$$

Thus we have shown

$$\lim u^{f(g_k,t_k)} g_k x_k y_k u^{-t_k} \in a^{-\log(1\pm\delta)}x^{[\frac{-\delta}{1-\delta},\frac{\delta}{1-\delta}]}$$

Now suppose that $|c_k| \leq |y_k|$. Thus $t_k=\frac{\delta}{|y_k|}$. That is $t_ky_k=\pm \delta$. As in the previous case, we wish to bound the quantity $\lim_{k \to \infty} d_k -t_kc_k$.  Because $|c_k|\leq |y_k|$, we have that 

$$ -\delta \leq t_k c_k \leq \delta$$ 

Thus, since $d_k \to 1$, we have that $\lim_{k \to \infty} [d_k - t_k c_k] \in [1 - \delta , 1 + \delta]$. 
So $\lim  \pi(u^{f(g_k,t_k)}g_kx_ky_ku^{-t_k})  $ is contained in 
$$a^{[-\log(1+\delta), -\log(1-\delta)]}$$ 

As before, $c_kx_k+d_ky_k \to 0$, and 
$a_k x_k + b_ky_k \to 0$, so

\begin{align*}\lim[(a_kx_k + b_ky_k)+f_k(c_kx_k+d_ky_k)] & =\lim[f_k(c_kx_k+d_ky_k)] \\
& = \lim[\frac{(t_ka_k-b_k)(c_kx_k+d_ky_k)}{d_k-t_kc_k}] \\
& = \frac{\lim[(t_ka_k-b_k)(c_kx_k+d_ky_k)]}{\lim[d_k-t_kc_k]}\end{align*}

The numerator equals
$$ \lim[t_ka_kc_kx_k + t_ka_kd_ky_k -b_kc_kx_k - b_kd_ky_k]$$

Because $x_k , y_k, b_k \to 0$, and $c_k, d_k \to 1$, the last two terms in the expansion of the numerator tend to $0$. 
We've already shown that $t_kc_k$ is bounded and since $a_k \to 1$ and $x_k \to 0$, the first term tends to $0$. 
Finally, $t_ky_k = \pm \delta$ and since $a_k, d_k \to 1$, we have that the numerator is $\pm \delta$. 

And as we've already shown, $\lim[d_k -t_kc_k] \in [1-\delta, 1+\delta]$. Altogether, we find that the limit of the $x$-coordinate of our vector lies in $$\Big[\frac{-\delta}{1-\delta}, \frac{-\delta}{1+\delta}\Big]\cup\Big[\frac{\delta}{1+\delta}, \frac{\delta}{1-\delta}\Big]$$   

This proves the claim. 

Now, we have that \begin{align*} u^{f(g_k,t_k)}M' & = u^{f(g_k,t_k)}g_kx_ky_k
M_k \\ & = \big[u^{f(g_k,t_k)}g_kx_ky_ku^{-t_k}\big]u^{t_k}M_k
\end{align*}

Therefore, $\lim_{k \to \infty} u^{f(g_k,t_k)}M' =\lim_{k \to
\infty} \beta_\delta u^{t_k}M_k$. From the compactness of $K$ it
follows that $\beta_\delta K \cap K \neq \emptyset$ and thus
$\beta_\delta \in \text{Stab}_{BX_0}(\mu)$ whenever \ref{e:1} and
\ref{e:2} are satisfied. Since $ \text{Stab}_{BX_0}(\mu)$ is
closed and \ref{e:1} and \ref{e:2} hold on arbitrarily large
subsets, we can let $\delta \to 0$ to see that $\beta_\delta \in
\text{Stab}_{BX_0}(\mu)$ and more precisely, that $\beta_\delta
\in \text{Stab}_{BX_0}(\mu)^\circ$.

\end{proof}

\section{Stabilizer is $U$: arclength measures}\label{S:st=u}

 For any $M \in \mathscr{L}_m$, let $M_X
\subseteq \mathscr{L}_m$ be the connected component of $M$ of the
space of all surfaces $xM \in \mathscr{L}_m$ for $x \in X$. If $M
\in \mathscr{L}_m$, we call $BM_X \subseteq \mathscr{L}_m$ the
\emph{$U$-normalizer space of }$M$.

Note that each $U$-normalizer space is $U$ invariant and that
$\mathscr{L}_m$ is a disjoint union of $U$-normalizer spaces.

This next lemma is essentially Proposition 1.6.10 of \cite{Mo}.

\begin{lemma}\label{l:u-norm} Suppose that
$\text{\emph{Stab}}_{BX_0}(\mu)^\circ = U$. Then $\mu $ is
supported on a $U$-normalizer space of a single surface.
\end{lemma}

\begin{proof}
If $\mu $ is supported on countably many $U$-normalizer spaces,
then it assigns positive measure to at least one. By ergodicity,
it would be supported on a single $U$-normalizer space. Thus, to
prove our claim we can assume that $\mu$ is supported on
uncountably many $U$-normalizer spaces, and then arrive at a
contradiction.

Assume that $\mu$ is supported on uncountably many $U$-normalizer
spaces. With $\Omega _\varepsilon$ as in Lemma~\ref{l:transverse},
$$\bigcup_{n \in \mathbb{N}} \Omega _{1/n}$$ is conull, so there
must be a fixed $n$ such that $\mu$ restricted to $\Omega _{1/n}$
is supported on uncountably many $U$-normalizer spaces.

Since uncountable sets contain a limit point, there is an $M' \in
\Omega _{1/n}$ and a sequence of $M_k \in \Omega _{1/n}$  such
that $M_k \to M'$ and $M_k \notin B(M')_X$. Since $M_k \to M'$,
there must exist a sequence of elements $(g_k,x_ky_k) \in
SL(2,\RR) \ltimes XY$ such that $(g_k, x_ky_k)M_k =M'$ with
$(g_k,x_ky_k) \to 1$.

Since $(g_k,x_ky_k) \notin BX$, we have for all $k$ that either $c_k
\neq 0$ or $y_k \neq 0$. Either way, we can pass to a subsequence
and apply Lemma~\ref{l:transverse} to find a nontrivial element of
$\text{Stab}_{BX_0}(\mu )^\circ-U$. This contradicts our
assumption that $\text{Stab}_{BX_0}(\mu )^\circ=U$.

\end{proof}

We will denote the identity component of the group of diagonal
matrices in $SL(2,\RR)$ by $A$.

\begin{lemma}\label{l:u-orbit} Suppose that
$\text{\emph{Stab}}_{BX_0}(\mu)^\circ = U$. Then there exists some
$a \in A$ and some $M' \in \mathscr{L}_m$ such that
$\text{\emph{Supp}}(\mu )=aUM'$.
\end{lemma}

\begin{proof}

Let $BM_X$ be the $U$-normalizer space from Lemma~\ref{l:u-norm}
and note that the proof of Lemma 3.3 from \cite{E-M-M} shows there
is a $U$-invariant, Borel subset ${\Omega _{AX}} \subset BM_X$,
such that
$$\mu ({\Omega _{AX}} )=1$$ and $${\Omega _{AX}} \cap g {\Omega _{AX}} =\emptyset \text{
for all } g \in AX_0 - \text{Stab}_{AX_0}(\mu )$$ Thus our claim
is that $\mu |_{\Omega _{AX}}$ is supported on $aUM'$ for a fixed
$a$ and $M'$. (Notice that $aUM'$ is a $U$-orbit since $a$
normalizes $U$, so $\text{Supp}(\mu | _{\Omega _{AX}} )$ certainly
contains some set of the form $aUM'$.)

Suppose $a_1 U M_1 \subseteq \text{Supp}(\mu | _{\Omega _{AX}} )$
for some $a _1 \in A$ and some $M_1 \in M_X$. We will show that
$a_1 U M_1 = \text{Supp}(\mu | _{\Omega _{AX}} )$ and thus will
prove our claim. But first it will be helpful to show that
$\text{Supp}(\mu | _{\Omega _{AX}} ) \subseteq BM_X$ can be
thought of as being arbitrarily narrow in the $M_X$-direction.

The surface $M_1$ is contained in $M_X$. Let $I^{M_1}_\varepsilon
\subseteq M_X$ be an interval of diameter $\varepsilon >0$ that
contains $M_1$. The set $BI^{M_1}_\varepsilon$
contains an open set in $BM_X$ that contains $a_1UM_1$. As the latter set is
contained in the support of $\mu$, we have that $\mu (
BI^{M_1}_\varepsilon )>0$. Because $A$ normalizes $U$, the set
$BI^{M_1}_\varepsilon$ is $U$-invariant. It follows from
ergodicity that $ BI^{M_1}_\varepsilon$, and thus ${\Omega _{AX}}
\cap BI^{M_1}_\varepsilon$, is conull.

Notice that $\varepsilon >0$ in the above paragraph was arbitrary.
We will assume that $\varepsilon$ is sufficiently small depending
on $a_1$. How it depends on $a_1$ will be explained below.
If $a_1 U M_1 \neq \text{Supp}(\mu | _{\Omega _{AX}} )$, then
there is a $U$-orbit contained in $\text{Supp}(\mu | _{\Omega
_{AX}}) \cap BI^{M_1}_\varepsilon$ that is distinct from
$a_1UM_1$. This orbit necessarily has the form $a_2 U M_2$ for
some $a_2 \in A$ and some $M_2 \in I^{M_1}_\varepsilon$. Thus, there is a $w \in X$ with norm less than $\varepsilon$ such that
$w M_1 =M_2$.

Recall that $A$ normalizes $X$, so $a_1wa_1^{-1} \in X$.
Furthermore, by choosing $\varepsilon > 0$ sufficiently small
depending on $a_1$, we may assume that $a_1wa_1^{-1} \in X_0$.

Let $$g=a_2 w a_1^{-1}$$ Note that $$g=a_2a_1^{-1}(a_1wa_1^{-1})
\in AX_0$$ Since $w$ commutes with $U$, we also have that
$$g(a_1UM_1)=a_2UM_2$$ The final piece of information we need
about $g$ is that $g \neq 1$. This follows from the fact that the
two $U$-orbits, $a_1UM_1$ and $a_2UM_2$, are distinct.

Altogether we have the following contradiction:
$$a_2UM_2 \subseteq g \text{Supp}(\mu | _{\Omega _{AX}} ) \cap
\text{Supp}( \mu | _{\Omega _{AX}} ) \subseteq g{\Omega _{AX}}
\cap {\Omega _{AX}} = \emptyset$$ The final equality in the line
above follows from the definition of ${\Omega _{AX}}$ and the fact
that $g \notin \text{Stab}_{AX_0}(\mu)^\circ =1$.

\end{proof}

\begin{lemma}\label{l:circles} Suppose that
$\text{\emph{Stab}}_{BX_0}(\mu)^\circ = U$. Then there is a
$U$-periodic surface $N \in \mathscr{L}_m$ such that $UN \subseteq
\mathscr{L}_m$ is homeomorphic to $S^1$ and $\mu$ is arc-length
measure on $UN$.
\end{lemma}

\begin{proof}
Let $N=aM'$ where $a$ and $M'$ are as in Lemma~\ref{l:u-orbit}.
Then $\mu$ is supported on the $U$-orbit of $N$. Because $U$
stabilizes $\mu$, the measure must descend from
Haar measure on $U$. Since the measure is a probability measure, the $U$-orbit is a closed circle, and $\mu$ is arc-length.

 \end{proof}

\section{Stabilizer is not unipotent: Entropy}\label{S:entropy}

By Lemma~\ref{l:circles}, we may now assume that
$\text{Stab}_{BX_0}(\mu)^\circ$ properly contains $U$. In this section, we assume there
is an element of $\text{Stab}_{BX_0}(\mu)^\circ$ that is not
unipotent.

There are three cases to be dealt with. First,
$\mu(\mathscr{L}_m^X)=1$ and $\mu$ is $X$-invariant. Second,
$\mu(\mathscr{L}_m^X)=1$ and $\mu$ is not $X$-invariant. Third,
$\mu(\mathscr{L}_m^X)=0$. But before we proceed with individual
cases, we will show that almost every surface in $\mathscr{L}_m$
admits relative translations of singularities by arbitrarily long
vertical vectors.

\subsection{Fibers of $\mathscr{L}_m$ and $\mu$.}\label{s:fibers}
Fix a surface $M \in \mathscr{L}_m$ and let $M_{XY} \subseteq
\mathscr{L}_m$ be the set of all surfaces in $\mathscr{L}_m$ that
can be realized as $v_1v_2 \cdot \cdot \cdot v_mM$ for some $v_1,
v_2, ...v_m \in XY$.

Let
$$\Gamma =\{\, g \in SL(2,\RR) \mid gM_{XY}=M_{XY}
\,\}$$

\begin{lemma}\label{l:fiber}
$\Gamma$ is discrete.
\end{lemma}

\begin{proof} If $\Gamma$ is not discrete then there
is a sequence $\{\gamma _n \} \subseteq \Gamma -1$ such that
$\gamma _n \to 1$. Thus, $\gamma _n M \to M$.

Since $\gamma _n M \in M_{XY}$, we have $\gamma _n M =
v_{1,n}v_{2,n} \cdot \cdot \cdot v_{m,n}M$ for some
$v_{1,n},v_{2,n},..., v_{m,n} \in XY$.

In a sufficiently small neighborhood of $M$, coordinates are given
for $\mathscr{L}_d$ by the absolute and relative homology. Since the absolute
homology of $v_{1,n}v_{2,n} \cdot \cdot \cdot v_{m,n}M$ and $M$
agree, the absolute homology of $\gamma _n M$ and $M$ agree which
implies that $\gamma _n M=M$. But the set of all $g \in
SL(2,\RR)$ such that $gM=M$ is discrete. Thus, the
sequence $\{\gamma _n \}$ is bounded away from $1$, a
contradiction.

\end{proof}

The space $\mathscr{L}_m$ fibers over
$SL(2,\RR)/\Gamma$ with fibers homeomorphic to
$M_{XY}$. By the previous lemma,
$SL(2,\RR)/\Gamma$ is a manifold so there is a
measure $\mu _{\pi}$ on $SL(2,\RR)/\Gamma$ and a
fiber measure $\mu _{M_{XY}}$ for every fiber $M_{XY}$ such that
$\mu$ is obtained by integrating the fiber measures over $\mu_\pi$.

\subsection{Vertical translations are conull}\label{s:vertconull}

For any line through the origin of the plane, $\ell \in
\mathbb{P}^1( \mathbb{R})$, we let
$$\mathscr{L}^\ell _m= \{\, N \in \mathscr{L}_m \mid  \ell N
\subseteq \mathscr{L}_m  \,\}$$

For any $g\Gamma \in SL(2,\RR)/\Gamma$, we define
$$\Sigma _{g \Gamma}^\ell=\{\,x\in gM_{XY} \mid x \notin
\mathscr{L}_m^\ell \,\}$$ so that $\mathscr{L}_m -
\mathscr{L}^\ell _m$ equals $$ \bigcup _{g\Gamma \in
G/\Gamma}\Sigma _{g \Gamma}^\ell$$

\begin{lemma}\label{L:count}

 Given $g \in SL(2,\RR)$, and $\ell, w \in \mathbb{P}^1(\RR)$, the intersection $\Sigma _{g \Gamma}^\ell \cap \Sigma _{g \Gamma}^w$ is uncountable only if $\ell =w$.

\end{lemma}

\begin{proof} Suppose $\{N_\alpha \}_{\alpha \in \mathscr{A}} \subseteq \Sigma _{g \Gamma}^\ell \cap \Sigma _{g \Gamma}^w$ where $\mathscr{A}$ is uncountable.

For any $\alpha \in \mathscr{A}$, let $x_\alpha \in \ell$ and $y_\alpha \in w$ be saddle connection on $N_\alpha$ between distinct singularities. Since $\mathscr{A}$ is uncountable, the set of triples $(N_\alpha,x_\alpha,y_\alpha) \in gM_{XY} \times \ell \times w$ contains an accumulation point $(N,x,y)$. We let $\{(N_i,x_i,y_i)\}\subseteq \{(N_\alpha,x_\alpha,y_\alpha)\}$ be a sequence that converges to $(N,x,y)$

Let $v_i \in \RR ^2$ be such that $v_i N_i =N$. Then $x -v_i$ and $x_i$ are saddle connections on $N_i$, that each converge to $x$. Thus, we may assume that $x-v_i = x_i$ for all $i$. Hence, $v_i \in \ell$ for all $i$. 
Similarly, $v_i \in w$ for all $i$, so $\ell =w$.
\end{proof}

The following lemma is essentially Lemma 5.4 from \cite{E-M-M}.

\begin{lemma}\label{L:vertconull}
 Let $\mu$ be an ergodic $U$-invariant measure on $\mathscr{L}_m$. Then
 $\mu(\mathscr{L}^Y_m)=1$
\end{lemma}

\begin{proof} Suppose $\mu (\mathscr{L}^Y_m)<1$. Then
by Fubini's theorem, there is a set $E\subseteq
SL(2,\RR)/\Gamma$ such that $\mu _{\pi}(E)>0$
and $\mu_{gM_{XY}}(\Sigma ^Y_{g\Gamma})>0$ for all $g \Gamma
\in E$.

Let $\lambda$ be Lebesgue measure on $U \cong \mathbb{R}$. By the
pointwise ergodic theorem there is a $g_0\Gamma \in
SL(2,\RR)/\Gamma$ and a set $U_0 \subseteq U$ such
that $\lambda (U_0)>0$ and $ug_0\Gamma \in E$ for all $u \in U_0$. Note that $\lambda(U_0)>0$ implies that $U_0$ is uncountable.

It can be checked that $u \Sigma ^{u^{-1}Y}_{g_0 \Gamma}= \Sigma
^Y_{ug_0\Gamma}$. This fact and the $U$-invariance of $\mu$
implies

$$\mu_{g_0M_{XY}}(\Sigma ^{u^{-1}Y}_{g_0 \Gamma})
= \mu_{ug_0M_{XY}}(u\Sigma ^{u^{-1}Y}_{g_0 \Gamma})
 = \mu_{ug_0M_{XY}}(\Sigma ^{Y}_{ug_0 \Gamma}) > 0$$
for all $u \in U_0$. In particular, if
$\widehat{\mu}=\mu_{g_0M_{XY}}$ and $\Sigma ^{\ell}_{g_0
\Gamma}=\Sigma ^\ell$, then there are uncountably many $\ell \in
\mathbb{P}^1(\mathbb{R})$ with $\widehat{\mu}(\Sigma ^\ell)
>0$.

We form a graph with a vertex for each $\ell \in
\mathbb{P}^1(\mathbb{R})$ with $\widehat{\mu}(\Sigma ^\ell)
>0$ and an edge for each distinct pair $\ell, w \in
\mathbb{P}^1(\mathbb{R})$ with $\widehat{\mu}(\Sigma ^\ell
\cap \Sigma ^w)>0$.

If a graph with uncountably many vertices has only countably many
edges, then there are uncountably many isolated vertices. Thus,
there must be uncountably many distinct pairs $\ell, w \in
\mathbb{P}^1(\mathbb{R})$ with $\widehat{\mu}(\Sigma ^\ell
\cap \Sigma ^w)
>0$, or else if $\mathcal{I}$ is the set of isolated vertices,
then $\widehat{\mu}(\cup _{\ell \in \mathcal{I}} \Sigma ^\ell
) =\infty$, which would be a contradiction.

If $\ell, w \in \mathbb{P}^1(\mathbb{R})$ are distinct and
$\widehat{\mu}(\Sigma ^\ell \cap \Sigma ^w)
>0$, then let $P(\ell,w)=\{\,N \in \Sigma ^\ell \cap \Sigma ^w
\mid \widehat{\mu}(\{N\})>0 \,\}$. Note that
$P(\ell,w)\neq\emptyset$ since $\Sigma ^\ell \cap \Sigma
^w $ is countable.

Form a second graph with a vertex for each $P(\ell,w)$ and edges
between $P(\ell,w)$ and $P(\ell ',w')$ if $P(\ell,w) \cap P(\ell
',w')\neq \emptyset$.

If $N \in g_0M_{XY}$, then there are only countably
many $\ell$ with $N\in \Sigma ^\ell$ since $N$ has a countable set
of saddle connections. Thus, only countably many $P(\ell, w)$
contain $N$. Furthermore, since the set of point masses for $\widehat{\mu}$ is countable, the edge set of our second graph is
countable, and therefore there are uncountably many pairs $\ell, w
\in \mathbb{P}^1(\mathbb{R})$ with $P(\ell,w)$ pairwise disjoint,
so the measure of the union of such sets is infinite. Thus, $\widehat{\mu}(g_0M_{XY})=\infty$. This is a
contradiction.
\end{proof}

\subsection{First case: Lebesgue}\label{s:leb}

We will assume in this section that $\mu (\mathscr{L}_m^X)=1$ and
that $\mu$ is $X$-invariant.

Since $UX$ is a codimension $1$ unipotent subgroup of the
non-unipotent group $\text{Stab}_{BX}(\mu)^\circ$, we have that
$BX = \text{Stab}_{BX}(\mu)^\circ$.
For $a^s \in A$, we let $h_{\mu}(a^s)$ be the entropy of the
transformation $a^s$ on $\mathscr{L}_m$ with respect to $\mu$.
Recall that $h_{\mu}(a^s)=h_{\mu}(a^{-s})$.

Since $\mathscr{L}_m$ is foliated by leafs that are locally the
orbits of $UX$, the entropy of $a^s$ is determined by the rate of
expansion in the $U$ and $X$ directions, or similarly, by the
expansion in the $U^t$ and $Y$ directions. Precisely, the proof of
Theorem 9.7 from \cite{M-T} yields

\begin{lemma}\label{l:entropy} Suppose $\mu$ is $BX$-invariant.
Then $h_{\mu}(a^s) = 3|s|$. Also, $h_\mu(a^{-s}) \leq 3|s|$ with
equality if and only if $\mu$ is $U^tY$-invariant.

\end{lemma}

From Lemma~\ref{l:entropy} we have

\begin{proposition}{\label{thm: vyinv}} If $\mu$ is $BX$-invariant, then it is  $U^tY$-invariant.
\end{proposition}

\begin{proof} We have $$3|s|=h_\mu (a^s)= h_\mu(a^{-s}) \leq
3|s|$$so the inequality is an equality.
\end{proof}

Now we have

\begin{proposition} If $\mu$ is $BX$-invariant,
 then $\mu$ is the unique ergodic,
 \slr-invariant measure with support equal to $\mathscr{L}_m$.
 \end{proposition}

\begin{proof}

By Proposition~\ref{thm: vyinv}, $\mu$ is also $U^t$-invariant, so
it is \slr-invariant since the subgroups $U$, $A$, and $U^t$
generate \slr.

McMullen classified the $SL(2,\RR)$-invariant
ergodic probability measures on the space of abelian differentials
in genus $2$; see Theorem 1.5 \cite{mc2}. It follows from the
classification that $\mu$ either equals $\nu _{10}$ or $\mu$ is as desired.
But $\nu _{10}$ only has a $3$-dimensional support, and $\mu$ has
support equal to $\mathscr{L}_m$ since the support of a measure is
closed and $\mu$ is invariant under $X$ and $Y$ as well as \slr.
\end{proof}

\subsection{Second case: decagon}\label{s:entropy-decagon}

We will assume in this section that $\mu (\mathscr{L}_m^X)=1$ and
that $\mu$ is not $X$-invariant.
Because $\text{Stab}_{BX}(\mu)^\circ$ is not unipotent, there is
some $x \in X$ such that $x^{-1}Ax \in
\text{Stab}_{BX}(\mu)^\circ$. Since $x$ commutes with $U$, the
measure $x_*\mu$ is ergodic and is invariant under $B$.

In order to apply an argument similar to the proof of
Proposition~\ref{thm: vyinv}, we first have to show that $x_*\mu$
can not detect the $X$-direction the $\mathscr{L}_m$. This will
affect the expansion of the $UX$-foliation under $a^s$ that is
visible to $\mu$, and thus will alter the calculation of
$h_{x_*\mu}(a^s)$.

\begin{lemma}\label{l:xblind} There is a conull $\Omega \subseteq
\mathscr{L}_m^X$ such that if $x_*\mu$ is not invariant under $X$,
then for all $M \in \Omega$, we have $$(UXM)\cap\Omega =
(UM)\cap\Omega$$ \end{lemma}

\begin{proof} Let $\Omega$ be as in Lemma 3.3 of \cite{E-M-M}, so that $\Omega$ is $U$-invariant and so that if $x \in X$ and $x\Omega \cap \Omega \neq \emptyset$ then $x = 1$. Thus if $u_1xM \in \Omega$ and $u_1xM=u_2M$ then $xM=u_1^{-1}u_2M \in \Omega$, so $x=1$.
\end{proof}

Analogously, we will need to know that the $Y$ direction does not contribute to the determination of $h_{x_*\mu}(a^{-s})$.

\begin{lemma}\label{l:yblind} There is a conull $\Omega \subseteq
\mathscr{L}_m^Y$ such that if $x_*\mu$ is not invariant under $X$,
then for all $M \in \Omega$, we have $$(U^tYM)\cap\Omega =
(U^tM)\cap\Omega$$ \end{lemma}

\begin{proof} The proof is essentially the proof of Proposition 5.5 from \cite{E-M-M}.

We reproduce it here for convenience.
Let $\Omega \subseteq \mathscr{L}_m^Y$ be a generic set such that
$x_*\mu(\Omega)=1$ and $a^sM \in \Omega_\varepsilon$ for most
$s>0$.
Let $vy \in U^tY$ and $M,M'\in \Omega$ be such that $vyM=M'$. We
wish to show that $y=0$. Suppose $y\neq 0$, and we will reach a
contradiction.
Choose a sequence of real numbers $s_k \to \infty$ such that
$a^{s_k}M,a^{s_k}M' \in \Omega _\varepsilon$ for all $k$. Note
that $a^{s_k}(vy)a^{-s_k} \to 0$ and that
$$||a^{s_k}va^{-s_k}|| \leq ||a^{s_k}ya^{-s_k}||$$ for $k \gg 0$.

Let $g_k=a^{s_k}va^{-s_k} \in SL(2,\RR)$ and $y_k
=a^{s_k}ya^{-s_k}\in Y$ and $t_k = \delta / ||a^{s_k}ya^{-s_k}||$.
By Lemma~\ref{l:transverse}, $
\text{Stab}_X(x_*\mu)^\circ \neq 0$. This is a contradiction.
\end{proof}

Using Lemmas~\ref{l:xblind} and~\ref{l:yblind}, one can apply the proof of Theorem 9.7
from \cite{M-T} to show

\begin{lemma}\label{l:entropyV} Suppose $x_*\mu$ is $B$-invariant and not $X$-invariant. Then
$h_{x_*\mu}(a^s) = 2|s|$. Also, $h_{x_*\mu}(a^s) \leq 2|s|$ with
equality if and only if $x_*\mu$ is $U^t$-invariant.
\end{lemma}

\begin{proposition}{\label{thm: vinv}} If
$x_*\mu$ is $B$-invariant, then it is $U^t$-invariant.
\end{proposition}

\begin{proof}
By Lemma~\ref{l:entropyV}, we have $$2|s|=h_{x_*\mu} (a^s)=
h_{x_*\mu}(a^{-s}) \leq 2|s|$$so the inequality is an equality.
\end{proof}

\begin{proposition} If $x_*\mu$ is $B$-invariant,
 then $x_*\mu=\nu _{10}$. \end{proposition}
\begin{proof}

By assumption, $x_*\mu$ is $U$-invariant and $A$-invariant. By
Proposition~\ref{thm: vinv}, $x_* \mu$ is also $U^t$-invariant.
Since the subgroups $U$, $A$, and $U^t$ generate \slr, $x_*\mu$ is
\slr-invariant. Thus, by Theorem 1.5 of \cite{mc2}, $x_*\mu$
equals $\nu _{10}$ or a measure whose support is $\mathscr{L}_m$.
By lemma~\ref{l:xblind}, $x_*\mu = \nu _{10}$.
\end{proof}

\subsection{Third case: Lattice surfaces on the boundary}\label{s:entropy-full-support}

For our final case, we assume that $\mu (\mathscr{L}^X_m)=0$. Thus
$\text{Supp}(\mu) \subseteq \text{HC}(r)$ for some $r>0$ from
which it follows that there exists some $w\in X$ with $|w|=r$ such that $w$ is a horizontal saddle connection for all $M \in \text{Supp}(\mu)$.

Since $\text{Stab}_{AX_0}(\mu)^\circ$ is not unipotent,  $\text{Stab}_{AX_0}(\mu)^\circ$ contains a neighborhood of the identity of a $1$-parameter subgroup of $AX$ that is not contained in $X$. Any such neighborhood is of the form $\{x^{-1}a^sx\} _{|s|<\varepsilon}$ for some $x\in X$.

Note that if $w$ is a horizontal saddle connection on $M$, then $e^sw-e^sx+x$ is a horizontal saddle connection on $x^{-1}a^sxM$.
Since $x^{-1}a^sx$ stabilizes $\mu$, it follows that $x^{-1}a^sx$ stabilizes the support of $\mu$ as well. Since saddle connections on a surface are discrete, it follows that $w=e^sw-e^sx+x$ and therefore $x=w$.

It can be easily checked that $x^{-1}U^tx$ acts on the space of surfaces in $\mathscr{L}_m$ that have $x$ as a horizontal saddle connection, so
$x^{-1}SL(2,\RR)x$ acts on $\mu$ with a stabilizer that includes $x^{-1}Bx$.

Similar to Lemma~\ref{l:entropyV}, $h_{\mu}(xa^{-s}x^{-1})=h_{\mu}(x^{-1}a^sx) = 2|s|$ which implies that $\mu$ is $x^{-1}SL(2,\RR)x$-invariant.

\section{Stabilizer is $UX$}\label{S:UX}

In this final section we will prove the following proposition which completes our proof of Theorem~\ref{T:main}.

\begin{proposition}{\label{thm: UX}} If $\mu$ is horizontally invariant, then $\mu$ equals the unique ergodic,
 \slr-invariant measure with support equal to $\mathscr{L}_m$.
\end{proposition}

The contents of this section were described to us by Alex Eskin and were motivated by Ratner's proof of Theorem 2 from \cite{r5}.

Let $\nu$ be the unique ergodic,
 \slr-invariant measure with support equal to $\mathscr{L}_m$. Let  $f: \mathscr{L}_m \rightarrow \mathbb{R}_{\geq 0}$ be a continuous, compactly supported function and let $\varepsilon >0$ be given. 
 We will prove that $$ \Big|\int_{\mathscr{L}_m} f \, d\mu - \int_{\mathscr{L}_m} f d\nu \Big| < \varepsilon$$ and thus prove Proposition~\ref{thm: UX}.

For any $T>0$ and any $M \in \mathscr{L}_m^X$ we let 
$$ \mathscr{A}_U(f,T)(M) = \frac{1}{T}\int_{0}^T f(u^t M) \, dt $$
$$ \mathscr{A}_{UX}(f,T)(M) = \frac{1}{T^\frac32} \int_{0}^{ \sqrt T} \int_0^T f(u^t x^s M) \, dt \, ds$$

\subsection{An ergodic theorem}\label{s:weakerg}

By the Birkhoff ergodic theorem, there is some $E' \subseteq \mathscr{L}_m$ with $\mu (E')=1$, and such that if $M \in E'$ then $\lim _{T \to \infty} \mathscr{A}_U(f,T)(M)=\int_{\mathscr{L}_m}f \, d\mu $.

Choose a sequence $\varepsilon _n \to 0$ such that $\Sigma _{i=1}^\infty\sqrt{\varepsilon _n}$ converges. 

For all $n$, there is some $E_n \subseteq \mathscr{L}_m$ with $\mu (E_n ) > 1 -\varepsilon _n$ and $T_n >0$ such that if $M\in E_n$ and $T \geq T_n$ then \begin{align} \Big| \mathscr{A}_U(f,T)(M)-\int_{\mathscr{L}_m}f \, d \mu \Big| < \varepsilon _n \end{align}

We let $$E'_n=\{M\in \mathscr{L}_m \mid ds(\{s \in [0,\sqrt{T_n}] \mid x^s M \in E_n\} ) \geq (1 -\sqrt{\varepsilon _n}) \sqrt{T_n} \}$$

\begin{lemma}{\label{l:dun}} If $M \in E'_n$ and $C =\max\{ \sup f \,,\, 1\}$, then $$\Big| \mathscr{A}_{UX}(f,T_n)(M) -\int _{\mathscr{L}_m}f \, d\mu \Big| < \sqrt{\varepsilon _n} (2C + 2)$$ \end{lemma}

\begin{proof}
Let $S_n=\{s \in [0,\sqrt{T_n}] \mid x^s M \in E_n\}$ and $S_n^c=\{s \in [0,\sqrt{T_n}] \mid x^s M \notin E_n\}$. Now note that both  $$\frac{1}{\sqrt{T_n}} \int _{s \in S_n} \mathscr{A}_{U}(f,T_n)(x^sM) \, ds$$ and $$(1- \sqrt{\varepsilon _n}) \int _{\mathscr{L}_m} f \, d\mu$$ lie in the interval bounded by $(1- \sqrt{\varepsilon _n}) \big( \int _{\mathscr{L}_m} f \, d\mu - \varepsilon _n \big)$ and  $\int _{\mathscr{L}_m} f \, d\mu + \varepsilon _n$. Thus,

\begin{align*} & \Big| \frac{1}{\sqrt{T_n}} \int _{s \in S_n} \mathscr{A}_{U}(f,T_n)(x^sM) \, ds - (1- \sqrt{\varepsilon _n}) \int _{\mathscr{L}_m} f \, d\mu \Big| \leq \\
&  \leq \Big| \int _{\mathscr{L}_m} f \, d\mu + \varepsilon _n - (1- \sqrt{\varepsilon _n}) \big( \int _{\mathscr{L}_m} f \, d\mu - \varepsilon _n \big) \Big| \\
&  \leq \sqrt{\varepsilon _n}  \int _{\mathscr{L}_m} f \, d\mu + \varepsilon _n (2- \sqrt{\varepsilon _n}) \\
&   < \sqrt{\varepsilon _n}(C+2) \end{align*}

Also note

$$\Big| \frac{1}{T_n^{\frac32}}\int_{s \in S_n^c} \int_0^{T_n}f(u^tx^sM) \, dt \, ds - \sqrt{\varepsilon _n} \int_{\mathscr{L}_m} f \, d\mu \Big| \leq \sqrt{\varepsilon _n} C$$

The lemma follows since $ \mathscr{A}_{UX}(f,T_n)(M)$ is the sum of  
$$  \frac{1}{\sqrt{T_n}} \int _{s \in S_n} \mathscr{A}_{U}(f,T_n)(x^sM) \, ds$$ and $$  \frac{1}{T_n^\frac32}\int_{s \in S_n^c} \int_0^{T_n}f(u^tx^sM) \, dt \, ds $$
\end{proof}

\begin{lemma}
\label{lemma:measure:En:prime}
If $F_n \subseteq \mathscr{L}_m$ is the complement of $E'_n$, then $\mu(F_n) \leq \sqrt{\varepsilon _n}$.
\end{lemma}

\begin{proof}
Let $\chi_n$ denote the characteristic function of the complement of
$E_n$. Since $\mu(E_n) > 1 - \varepsilon_n$, we have
$$\int_{\mathscr{L}_m} \chi_n (M) \, d\mu (M) \leq \varepsilon_n$$
Since $\mu$ is invariant under $X$, we have for all $0 \le s \le \sqrt{T_n}$
that
$$\int_{\mathscr{L}_m} \chi_n(x^s M) \, d\mu(M) \leq \varepsilon_n$$
Integrating the inequality above with respect to $s$ gives
\begin{align*}\int_{\mathscr{L}_m} \left( \int_0^{\sqrt{T_n}} \chi_n(x^s M) \, ds \right)\, d\mu(M) &= 
\int_0^{\sqrt{T_n}} \left( \int_{\mathscr{L}_m} \chi_n(x^s M) \, d\mu(M) \right)\, ds \\
& \leq
\varepsilon_n \sqrt{T_n} \end{align*} 

Notice that $M \not\in E_n'$ exactly when
$$ \int_0^{\sqrt{T_n}} \chi_n(x^s M) \, ds \geq \sqrt{\varepsilon _n}
\sqrt{T_n}$$ Thus, the lemma follows from the third inequality of this proof.
\end{proof}

\begin{proposition}
\label{prop:weak:ergodic:theorem}
There exists a sequence $T_n \to \infty$ and 
\begin{itemize}
\item[{\rm (i)}] a subset $E_\mu \subset \mathscr{L}_m$
with $\mu(E_\mu) = 1$ such that for $M \in E_\mu$, 
\begin{equation*}
\lim_{n \to \infty} \mathscr{A}_{UX}(f,T_n)(M) = \int_{\mathscr{L}_m} f \, d\mu
\end{equation*}
\item[{\rm (ii)}] a subset $E_\nu \subset \mathscr{L}_m$
with $\nu(E_\nu) = 1$ such that for $M \in E_\nu$, 
\begin{equation*}
\lim_{n \to \infty} \mathscr{A}_{UX}(f,T_n)(M) = \int_{\mathscr{L}_m} f \, d\nu
\end{equation*}
\end{itemize}
\end{proposition}
 
\begin{proof}  By Lemma~\ref{lemma:measure:En:prime}, 
$\sum \mu(F_n)$
converges.
Let $E_\mu$ be the set of all $M\in\mathscr{L}_m$ such that $M$ is contained in at most finitely many of the $F_n$. By the 
Borel-Cantelli lemma, $\mu (E_\mu)=1$. Also note that for any $M \in E_\mu$, there is some $k$ such that $M \in E_n'$ whenever $n \geq k$. Thus, (i) follows from Lemma~\ref{l:dun}.

To prove (ii) note that the only properties of $\mu$ used were
ergodicity with respect to $U$ and $X$-invariance. These properties
are shared by $\nu$ as well. To ensure that the
sequence $T_n$ is the same for both $\mu$ and $\nu$ we choose the
$T_n$ so that -- in addition to being large enough to satisfy all properties above -- there also exists a
set $E''_n$ with $\nu(E''_n) > 1 - \varepsilon_n$, such that
$$
\Big| \mathscr{A}_U(f,T) - \int_{\mathscr{L}_m} f \, d\nu \Big| < \varepsilon_n$$
for all $T \geq T_n$ and all $x \in E''_n$. The rest of the proof of
(ii) is identical to that of (i). 
\end{proof}

\begin{corollary}
\label{cor:uniform:convergence:to:different:averages}
Let $f$ and $\varepsilon$ be as in the beginning of this section, and let
 $T_n$ be as in
Proposition~\ref{prop:weak:ergodic:theorem}.  Then for every $\delta'
> 0$ there
exist subsets
$\hat{E} \subset \mathscr{L}_m$ and $\tilde{E} \subset \mathscr{L}_m$ and an integer $n_0$ such that
\begin{itemize}
\item[{\rm (i)}] $\mu(\hat{E})> 1 - \delta'$ and $\nu(\tilde{E}) > 1 -
  \delta'$.  
\item[{\rm (ii)}] For $n > n_0$, and $M \in \hat{E}$, $|\mathscr{A}_{UX}(f,T_n)(M)
  - \int_{\mathscr{L}_m} f \, d\mu | < \varepsilon/4$. 
\item[{\rm (iii)}] For $n > n_0$, and $N \in \tilde{E}$,
  $|\mathscr{A}_{UX}(f,T_n)(N) - \int_{\mathscr{L}_m} f \, d\nu | < \varepsilon/4$. 
\end{itemize}
\end{corollary}

\subsection{Recurrence to compact sets}\label{s:recurrence}

The below theorem follows directly from Theorem H2 of Minsky-Weiss' \cite{M-W}.
\begin{theorem}
\label{theorem:returning:to:fixed:compact:set}
For any $\delta > 0$ there exists a compact set $K \subset \mathscr{L}_m$ such
that if $\lambda$ is any $U$-invariant probability measure on $\mathscr{L}_m$ that assigns measure $0$ to the set of surfaces that contain a horizontal saddle connection, then $\lambda(K) > 1 - \delta$. 
\end{theorem}

Recall that $a^t = \begin{pmatrix} e^t & 0 \\ 0 & e^{-t} \end{pmatrix}$. The proof of the next proposition is credited to  Elon
Lindenstrauss and Maryam Mirzakhani.

\begin{proposition}
\label{prop:support:recurrent:points}
Let $\delta > 0$ and let $t_n \to \infty$ be any sequence. Suppose $K$ is as in the previous theorem and that $\lambda$ is a $U$-invariant probability measure on $\mathscr{L}_m$  that assigns measure $0$ to the set of surfaces that contain a horizontal saddle connection. Let $F$ denote the set of $M
\in K$ such that there exists a subsequence $\tau_n$ of $t_n$,
with $a^{-\tau_n} M \in K$ for all $n$. Then $\lambda(F) > 1-\delta$. 
\end{proposition}

\begin{proof} Because $A$ normalizes $U$, we can apply
Theorem~\ref{theorem:returning:to:fixed:compact:set} to each pushforward measure $(a^{-t_n})_* \lambda$. Thus, for all $n \in \mathbb{N}$, 
$$\lambda( a^{t_n} K ) > 1 - \delta$$

 Therefore,
$$\lambda\left( \bigcap_{n=1}^\infty \bigcup_{k=n}^\infty a^{t_k} K \right)
\geq 1 - \delta$$ The lemma follows. 
\end{proof}

\subsection{An adaptation of an argument of Ratner's}\label{s:ratnerargument}

For any $M\in \mathscr{L}_m$, there is some $D>0$ such that $$a^{[-D,D]}y^{[-D,D]}v^{[-D,D]}M \subseteq \mathscr{L}_m$$ We let $W^D(M)=a^{[-D,D]}y^{[-D,D]}v^{[-D,D]}M$.

Choose $\delta > 0$ with $\delta \ll 1$, and let $K$ and $F$ be as in Theorem~\ref{theorem:returning:to:fixed:compact:set} and Proposition~\ref{prop:support:recurrent:points}. Since $K$ is compact, we may choose $D>0$ uniformly so that the above inclusion is satisfied for all $M \in K$. Furthermore, with $\varepsilon$ and the compactly supported $f$ as in the beginning of this section, we may assume that $D\ll 1$ is such that for any $M \in \mathscr{L}_m$, if there is some $N =a^{t_1}y^{t_2}v^{t_3}M$ where $|  t_i | \leq D$, then \begin{align}\label{eq:uc}
| f(M)-f(N) | < \varepsilon /8
\end{align}

There are constants $R \ll 1$ and $D'\ll D$ uniform for all $M \in K$ such that $W^{D'}(u^tx^sM)$ is a well-defined subset of $\mathscr{L}_m$ for $0\leq t,s \leq \sqrt{R}$, and such that if $N \in W^{D'}(M)$ and $0\leq t,s \leq \sqrt{R}$, then $u^{\gamma(N,t,s)}x^{\beta(N,t,s)}N\in W^D(u^tx^sM)$ for some $\gamma(N,t,s)\,,\,\beta(N,t,s)>0$. Note that the map $\phi : [0,\sqrt{R}] ^2 \rightarrow (\mathbb{R}_{\geq 0})^2$ given by $\phi (t,s) = \big(\gamma(N,t,s)\,,\,\beta(N,t,s) \big)$ is a diffeomorphism of $[0,\sqrt{R}]^2 $ onto its image, and for any $\eta >0$ we can choose $D' \ll 1$ such that 
\begin{align}\label{eq:vol} \Big|  \frac{\text{volume}_{UX}(J)}{\text{volume}_{UX}(\phi(J))} -1 \Big| < \eta
 \end{align}
for any Borel set $J \subseteq [0,1] ^2$.

For $M \in K$, we let \begin{align*} V^{D'} (M) & =\{\, u^{\gamma(N,t,s)}x^{\beta(N,t,s)}N \mid  N \in W^{D'} (M) \text{ and } 0 \leq s, t \leq \sqrt{R}\varepsilon /C \,\}  \\
 & \subseteq \bigcup_{0 \leq s, t \leq \sqrt R}W^D (u^tx^sM)
\end{align*}
By the compactness of $K$ there is some $\delta _0>0$ such that $\nu \big( V^{D'} (M) \big)>\delta_0$ for all $M \in K$.

Choose $\delta ' < \min \{ \delta , \delta _0 \}$, let $\hat{E}$ and $\tilde{E}$ be as in Corollary~\ref{cor:uniform:convergence:to:different:averages}, and let $M \in F \cap \hat{E}$.

Let $T_n$ be as in Proposition~\ref{prop:weak:ergodic:theorem} and let $t_n = \log (T_n/R)/2$.  We let $\tau _n$ be the subsequence of $t_n$ associated to $M$ as guaranteed  by Proposition~\ref{prop:support:recurrent:points}. We let $W^{D'}_{\tau_n}(M)=a^{\tau_n}W^{D'}(a^{-\tau_n}M)$ and $V^{D'}_{\tau_n}(M)=a^{\tau_n}V^{D'}(a^{-\tau_n}M)$.  Notice that the $A$-invariance of $\nu$ implies  \begin{equation}\label{eq:d0} \nu \big( V^{D'}_{\tau_n}(M) \big)>\delta_0 \end{equation} for all $\tau_n$.

If $N \in W^{D'}_{\tau_n}(M)$ --- so that $a^{-\tau_n}N \in W^{D'}(a^{-\tau_n}M)$ --- and $0\leq t' \leq T_n$ and $0 \leq s' \leq \sqrt{T_n}$, then we define $\Gamma(N,t',s')\,,\,B(N,t',s')>0$ as $$\Gamma(N,t',s')=\gamma(N,t'R/T_n,s'\sqrt{R}/\sqrt{T_n})T_n/R$$ and $$B(N,t',s')=\beta(N,t'R/T_n,s'\sqrt{R}/\sqrt{T_n})\sqrt{T_n}/\sqrt{R}$$ From \ref{eq:vol} it follows that the map $\phi _n: [0,T_n] \times [0,\sqrt{T_n}] \rightarrow (\mathbb{R}_{\geq 0})^2$ given by $\phi _n(t',s') = \big(\Gamma(N,t',s')\,,\,B(N,t',s') \big)$ satisfies 
\begin{align}\label{eq:vol2} \Big|  \frac{\text{volume}_{UX}(J)}{\text{volume}_{UX}(\phi_n(J))} -1 \Big| < 
\eta \end{align}
for any Borel set $J \subseteq [0,T_n] \times [0,\sqrt{T_n}] $.

By \ref{eq:d0} and our choice of $\delta '$, there is an $L \in V^{D'}_{\tau_n}(M)\cap \tilde{E}$.
Notice that $L=a^{\tau_n}u^{\gamma (a^{-\tau_n}N,t,s)}x^{\beta (a^{-\tau_n}N,t,s)}a^{-\tau_n}N$ for some $N\in W_{\tau_n}^{D'} (M)$.

 We have \begin{align*} u&^{\Gamma(N,t',s')}x^{B(N,t',s')}N   = \\
& = u^{\gamma(N,t'R/T_n,s'\sqrt{R}/\sqrt{T_n})T_n/R}x^{\beta(N,t'R/T_n,s'\sqrt{R}/\sqrt{T_n})\sqrt{T_n}/\sqrt{R}}a^{\tau_n}a^{-\tau_n}N \\
 & = a^{\tau_n} u^{\gamma(N,t'R/T_n,s'\sqrt{R}/\sqrt{T_n})} x^{\beta(N,t'R/T_n,s'\sqrt{R}/\sqrt{T_n})}a^{-\tau_n}N \\ 
 & \in  a^{\tau_n} W^D( u^{t'R/T_n}x^{s'\sqrt{R}/\sqrt{T_n}}a^{-\tau_n}M) \\ 
 & = a^{\tau_n} W^D( a^{-\tau_n}u^{t'}x^{s'}M) \\ 
  & = W_{\tau_n} ^D( u^{t'}x^{s'}M) \\
 & = a^{[-D,D]}y^{[-D\sqrt{R}/\sqrt{T_n},D\sqrt{R}/\sqrt{T_n}]}v^{[-DR/T_n,DR/T_n]}u^{t'}x^{s'}M \\
 & \subseteq W^D(u^{t'}x^{s'}M) \end{align*}

It now follows from \ref{eq:uc} that
\begin{equation*}
\label{eq:close:flow}
\Big|\mathscr{A}_{UX}(f,T_n)(M) - \frac 1{T_n^{\frac 32}} \int_0 ^{\sqrt{T_n}} \int_0^{T_n} f(u^{\Gamma(N,t',s')}x^{B(N,t',s')}N ) dt'ds' \Big|<  \varepsilon /8
\end{equation*}

and then from \ref{eq:vol2} that for $n \gg 0$

\begin{equation*}
\Big|\mathscr{A}_{UX}(f,T_n)(N) - \frac 1{T_n^{\frac 32}} \int_0 ^{\sqrt{T_n}} \int_0^{T_n} f(u^{\Gamma(N,t',s')}x^{B(N,t',s')}N ) dt'ds' \Big|< \varepsilon /8 
\end{equation*}

The two above inequalities give

\begin{equation}\label{eq:done}
\Big|\mathscr{A}_{UX}(f,T_n)(M) - \mathscr{A}_{UX}(f,T_n)(N)\Big| < \varepsilon /4
\end{equation}

Also note that for $n \gg 1$ and with $R \ll 1$
\begin{equation}
\label{eq:flowbox:same:average}
|\mathscr{A}_{UX}(f, T_n)(N) - \mathscr{A}_{UX}(f, T_n)(L)| < \epsilon/4
\end{equation}

Recall that $M \in \hat{E}$ and that $L \in \tilde{E}$. In view of   \ref{eq:done}, \ref{eq:flowbox:same:average}, and of
Corollary~\ref{cor:uniform:convergence:to:different:averages} (ii) and
(iii), we have  $$ \Big|\int_{\mathscr{L}_m} f \, d\mu - \int_{\mathscr{L}_m} f d\nu \Big| < \varepsilon$$  

This completes the proof of Proposition~\ref{thm: UX}.

\end{document}